\documentclass[12pt]{amsart}
\usepackage{amssymb,amscd,amsthm,verbatim,nicefrac,amsmath,color,fancyhdr,mathrsfs,braket,graphicx, turnstile,hyperref,indentfirst,csquotes,url,enumitem,amsfonts, mathtools,lipsum}
\usepackage[letterpaper, left=2.5cm, right=2.5cm, top=2.5cm,
bottom=2.5cm,dvips]{geometry}
\hypersetup{
    colorlinks=true,
    linkcolor=blue,
    filecolor=magenta,      
    urlcolor=cyan,
    citecolor=blue,
    pdftitle={Minimising the Peak to $p$-average ratio in the vectorial case},
    pdfpagemode=FullScreen,
    }

\usepackage{amsfonts}
\usepackage{amsmath}
\usepackage{amsthm}
\usepackage{amssymb}
\usepackage[english]{babel}
\usepackage{graphics}
\usepackage{latexsym}
\usepackage{longtable} \usepackage{mathrsfs}
\usepackage{graphicx}
\usepackage{wrapfig} 

\newtheorem{theorem}{Theorem}
\newtheorem{corollary}[theorem]{Corollary}
\newtheorem{lemma}[theorem]{Lemma}

\newtheorem{claim}[theorem]{Claim}
\newtheorem{example}[theorem]{Example}

\theoremstyle{definition}
\newtheorem{definition}[theorem]{Definition}
\newtheorem{remark}[theorem]{Remark}

\newcommand{\mA}{\mathcal{A}}

\newcommand{\mF}{\mathcal{F}}

\newcommand{\mP}{\mathcal{P}}
\newcommand{\mM}{\mathcal{M}}
\newcommand{\mO}{\mathcal{O}}
\renewcommand{\S}{\mathcal{S}}
\newcommand{\mY}{\mathscr{Y}}
\newcommand{\mD}{\mathcal{D}}

\newcommand{\mK}{\mathbb{K}}

\newcommand{\A}{\mathrm{A}}
\newcommand{\R}{\mathbb{R}}
\newcommand{\N}{\mathbb{N}}
\newcommand{\mS}{\mathbb{S}}

\newcommand{\X}{\textbf{X}}

\renewcommand{\P}{\mathrm{P}}

\newcommand{\D}{D}

\newcommand{\noi}{\noindent}

\newcommand{\al}{\alpha}
\newcommand{\be}{\beta}

\newcommand{\de}{\delta}
\newcommand{\De}{\Delta}
\newcommand{\e}{\varepsilon}

\newcommand{\la}{\lambda}

\newcommand{\Om}{\Omega}

\newcommand{\av}{ {\,\,-\hspace{-13pt}\int\!\!} }

\newcommand{\weak }{ \hspace{1pt}-\!\!\!\rightharpoonup}
\newcommand{\weakstar }{ \overset{\, *_{\phantom{|}}}{{\smash{\weak }}\, } }

\newcommand{\larrow}{\longrightarrow}
\newcommand{\ot}{\otimes}

\newcommand{\ri}{\rightarrow}

\newcommand{\p}{\partial}
\newcommand{\sub}{\subseteq}
\renewcommand{\set}{\setminus}
\newcommand{\by}{\times}
\newcommand{\rk}{\mathrm{rk}}

\newcommand{\Lip}{\mathrm{Lip}}

\newcommand{\ess}{\mathrm{ess}}

\newcommand{\Div}{\mathrm{Div}}

\newcommand{\bt}{\begin{theorem}}\newcommand{\et}{\end{theorem}}
\newcommand{\bd}{\begin{definition}}\newcommand{\ed}{\end{definition}}
\newcommand{\bl}{\begin{lemma}}\newcommand{\el}{\end{lemma}}
\newcommand{\beq}{\begin{equation}}\newcommand{\eeq}{\end{equation}}
\newcommand{\bc}{\begin{claim}}\newcommand{\ec}{\end{claim}}
\newcommand{\bex}{\begin{example}}\newcommand{\eex}{\end{example}}
\newcommand{\bcor}{\begin{corollary}}\newcommand{\ecor}{\end{corollary}}
\newcommand{\bp}{\begin{proof}}\newcommand{\ep}{\end{proof}}

\numberwithin{equation}{section}
\begin{document}

\title[Vectorial Calculus of Variations in $L^\infty$ and Aronsson systems]{An introduction to the vectorial Calculus of Variations in $\textbf{L}^\infty$ AND Aronsson PDE systems}

\author{Hussien Abugirda and Nikos Katzourakis}

\address{H.A.\ Department of Physics, College of Education for Pure Sciences, University of Kerbala, Kerbala, Iraq}
\email{hussien.a@uokerbala.edu.iq}


\address[corresponding author]{N.K.\ Department of Mathematics and Statistics, University of Reading, Whiteknights Campus, Lepper Lane, Reading RG6 6AX, UK}
\email{n.katzourakis@reading.ac.uk}


\date{16th of August 2026}


\thanks{\!\!\!\!\!\!\!\texttt{N.K.\ has been partially financially supported through the EPSRC grant EP/X017109/1.}}

\thanks{Invited paper for the conference proceedings of the \emph{Panhellenic Conference in Mathematical Analysis PCMA 2025}, held in Athens (Greece) 18-20 December 2025}

\begin{abstract} In this expository article, which is partially based on the lecture notes of a short course delivered by the second appearing author, we introduce the subfield of the Calculus of Variations that is concerned with the study of vectorial variational problems for supremal functionals, with emphasis on the associated PDE systems arising as extremality conditions. The scalar theory has a rather short history, first arising in the work of G.\ Aronsson in the 1960s. The vectorial case is even more recent and first arose in work of the second appearing author in the early 2010s. The Calculus of Variations in $L^\infty$ is an all-important field for numerous applications, but it presents serious difficulties and requires new tools, different to those required for the classical case of integral functionals and Euler-Lagrange equations.
\end{abstract}
\maketitle

\tableofcontents

\section{Introduction} \label{section1}

The \emph{Calculus of Variations} is one of the most prevalent and ubiquitous fields of mathematical analysis, having its roots in the work of Euler in the 1600s. In the 20th century, its rigorous development has gone hand-in-hand with associated breakthroughs in functional analysis, measure theory and nonlinear partial differential equations. 

Conventionally, the central object of study of the Calculus of Variations is functionals defined as integrals which act on an admissible class of functions defined on a domain, and typically forming a vector space. The integral functional usually depends through an energy density on the admissible functions themselves, and perhaps any number of its derivatives and the points of the domain. One then cares about studying variational problems involving these functionals, usually minimising or maximising them over an admissible class of functions, which includes also boundary or other side conditions or constraints (see e.g.\ \cite{D, FL}). 

It is rather remarkable there is a natural system of nonlinear Partial Differential Equations (PDEs) which arises in connection to variational problems for integral functionals, known as the \emph{Euler-Lagrange PDEs}. The ``quintessence of the Calculus of Variations" is the study of variational problems as well as their connection to the solvability of the Euler-Lagrange PDEs. Despite having classical roots in the Renaissance, this is still one of the most active fields of research nowadays (see e.g.\ \cite{D, Ev}). 

In the 1960s (some 300 plus years later!), the Swedish mathematician G.\ Aronsson was the first to consider variational problems for functionals defined as the (essential) supremum of an energy density, rather than an integral \cite{A1}-\cite{A7}. These objects are called \emph{supremal (or $L^\infty$) functionals}, and there and then was the birth of the field called nowadays \emph{Calculus of Variations in $L^\infty$} (see e.g.\ \cite{B, BL, C, K7} for the development of the scalar case). Despite the obvious similarity to integral functionals, this is more or less where shared properties stop. The Calculus of Variations in $L^\infty$ presents unprecedented difficulties, so crucial that essentially renders it a separate field. Practically, no  integral method or technique works for supremal functionals. Even the standard notion of minimisers is not appropriate in the $L^\infty$ setting and needs to be modified, whilst the classical theories of weak and distributional solutions to PDE do not apply either! 

Even though the Calculus of Variations in $L^\infty$ is a non-classical field, partly because supremal functionals do not arise in connection to conservation or least actions principles in physics, they are extremely important for a variety of applications inside and outside mathematics. Aronsson motivation was the problem of optimal Lipschitz extensions, whilst there are also connections to the problem of optimising quasi-conformal mappings \cite{CR, K5}, as well as curvature minimisation problems in geometry, relating to the Yamabe and Nirenberg problems \cite{KM1, KM3, MS, S}. More generally, it is crucial in numerous applications for which minimisation of the maximum energy e.g.\ modelling an error or misfit of a problem, is more realistic to the minimisation of the average, as spikes or large pointwise energy deviation are at the outset excluded. Example include PDE-constrained optimisation, optimal shape design, inverse problems, data assimilation, electromobility, engineering automotive applications, image processing, fracture, chemistry, and more (\cite{AK, BEJ, CKM, CK, ETT, GNP, K11, K12, K13, K14}).

In this expository paper we present the milestones of the development of the vectorial case of the Calculus of Variations in $L^\infty$. Vectorial problems in $L^\infty$ have a much shorter history, first arising in the work of the same author fifteen years ago in the early 2010s \cite{K1}. Even though the scalar case is nowadays well-developed and largely understood, the vectorial case presents several additional layer of complications and not even state-of-the-art scalar-valued $L^\infty$ techniques work. For this reason, we focus the exposition mostly on the properties of classical solutions to Aronsson-type PDE systems, but we discuss also variational characterisations and appropriate theories of generalised solutions. For first order supremal functionals (involving the admissible map plus derivatives up to first order), the Aronsson system is a second order quasilinear degenerate system PDE system with discontinuous coefficients, and in general it is not (degenerate) elliptic. 

The Aronsson system arises in connection to extremality conditions for such functionals. Even its very emergence is a curious and non-trivial matter, as supremal functionals are not differentiable, and therefore one cannot simply derive it by ``taking variations". Further, it is rather mysterious how it is possible for non-differentiable functionals to have necessary PDE extremality conditions, which in special cases are also be sufficient for minimality! The archetypal and simplest example of Aronsson system is the $\infty$-Laplacian, associated with the $L^\infty$ norm of the gradient. The scalar case of $\infty$-harmonic functions has attracted immense mathematical interest in the last 30 years, offering one of the most significant applications of the theory of viscosity solutions to second order PDE. The $\infty$-Laplacian is the formal limit case of the $p$-Laplacian as $p$ tends to infinity, but the breath and depth of the $L^\infty$ theory can only be compared to that of the usual ($2$-) Laplacian.

The present notes, which are partially based on a short course the second appearing author delivered at BCAM (Bilbao, Spain) in 2017-18, are intended as an introduction to (post-) graduate students of the mathematical sciences, as well as to non-experts in the field. They only require some basic mathematical maturity and some general knowledge of analysis, even though some more advanced parts, for instance about generalised solutions, do require some more familiarity with functional analysis, measure theory and nonlinear PDEs. The content consists of two main parts. After this introduction, in Section \ref{section2} we present some rudiments from integral calculus of variations, to illustrate the main ideas. In Section \ref{section3} we expound on the main theme, which is the theory of Aronsson PDE systems, as well as their connection to $L^\infty$ variational problems. There are no proofs  included, but instead we provide references where the relevant results have appeared. We will not discuss the higher order case of Aronsson PDE systems either, for which we refer for some explorative work to \cite{DK, KP}.


\section{Some elements from the integral Calculus of Variations} \label{section2}

\subsection{The $\infty$-dimensional sibling of optimisation}\label{subsection2.1} Let $\X$ be a finite-dimensional vector space, let also $E: \X \longrightarrow \R$ be a real valued differentiable functional (see e.g.\ Figure 1). 
\[
\underset{\scriptstyle{\text{Figure 1.  }}}{ \includegraphics[scale=.35]{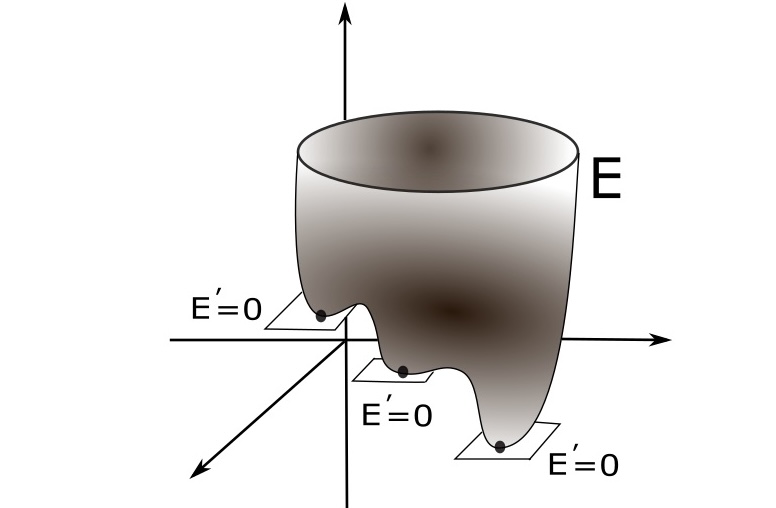}}
\]
If $E$ has a global minimum  at a point $x_0 \in \X$, namely if
\[
E (x_0)= \underset{x \in \R^n}{\inf} E(x), 
\]
then we have $\big(E(x_0+ te)- E(x_0)\big)/t\geq 0$ for all $t> 0$ and $e\in\X$, from where by taking the limit we infer that
\[
E'(x_0)=0.
\]
The same can be inferred if $x_0$ is merely a local extremum (i.e. maximum or minimum) at a point $x_0 \in \X$, namely either $E(x)\geq E(x_0)$ or $E(x)\leq E(x_0)$ for all $x$ in a ball centred at $x_0$. If in addition $E$ is twice differentiable, then we infer that
\[
E''(x_0) \geq 0, 
\]
if the point $x_0$ is a minimum, or $E''(x_0) \leq 0$ if the point $x_0$ is a maximum. In the Calculus of Variations, we aim at studying the extrema of functionals like $E$, but in infinite dimensions.

\subsubsection{Example}\label{subsection2.1.2} Consider the vector space $\X= W^{1,2}(\Om)$ (or $C^2(\overline{\Om})$), where $\Om \sub \R^n$ is open and bounded. Consider the problem of minimising the functional $E:  W^{1,2}(\Om) \longrightarrow \R$ defined as:
\[
E(u)=\int_\Om \Big(\frac{1}{2}|D u|^2 - f \Big), \ \ \ u\in W^{1,2}(\Om),
\]
where $D u$ symbolises the gradient $Du=(D_1u,...D_nu)$, $\D_i \equiv \p /\p x_i$, and suppose also $f\in L^2(\Om)$ (or $C(\overline{\Om})$). Then, by arguing as above, any minimiser $u$ of $E$ are the solutions to the Poisson equation (see e.g.\ \cite{Ev, D} and the next subsection for details):
\[
\ \ \ \De u=f , \ \ \ \text{in}  \  \Om.
\] 
\subsubsection{Example}\label{subsection2.1.3} Consider the space $\X= W^{1,p}(\Om; \R^N)$, where $\Om \sub \R^n$ is open and $p \in (1,\infty)$. More generally, consider the problem of minimising $E:  W^{1,p}(\Om; \R^N) \longrightarrow \R$ defined as:
\beq
\label{2.0}
E(u)=\int_\Om H(x,u(x),D u(x))\, dx, \ \ \ u\in W^{1,p}(\Om; \R^N).
\eeq
 Here $H: \Om \by \R^N \by \R^{N \by n} \longrightarrow \R$ is applied to Sobolev mappings  $u:\R^n \supseteq \Om\longrightarrow \R^N$, $u=(u_1,u_2,...,u_N)^\top$ and $D u$ symbolises the gradient matrix (also called ``Jacobian", but not by PDE theorists (!)):
\[
\ \ \ \ D u(x) = \big(\D_i u_\al(x)\big)_{i=1...n}^{\al=1...N} \, \in\, \R^{N\by n}.
\]
The appropriate $p$ depends on the assumptions. For $p\in[2,\infty)$, a special case is the celebrated $p$-Laplacian is the divergence system
\beq \label{1.33}
\ \ \ \De_p u:= \Div \big(|D u|^{p-2}D u \big) = 0\ \ \text{ in }\Om,
\eeq
and comprises the Euler-Lagrange equation which describes extrema of the model $p$-Dirichlet integral functional
\beq  \label{1.34}
\ \ \ E_p(u) := \int_\Om|D u|^p, \ \ \ u\in W^{1,p}(\Om;\R^N),
\eeq
Above and subsequently, for any $X\in \R^{N\by n}$, the notation $|X|$ symbolises its Euclidean (Frobenius) norm: 
\[
|X| = \left(\sum_{\al=1}^N\sum_{i=1}^n\, (\X_{\al i})^2\right)^{\!\! 1/2},
\]
and the associated inner product is given by $X:Y :=\mathrm{tr}(X^\top Y)$, for any $X ,Y \in \R^{N \by n}$. The pair \eqref{1.33}-\eqref{1.34} is of paramount important and has been studied extensively in the literature, due to both its intrinsic interest and also its numerous applications.

\subsection{Critical points and Euler-Lagrange equations }\label{subsection2.2} Consider the functional $E$ from \eqref{2.0}, and let $u \in C^2(\Om;\R^N)$ be a vector function defined on some open bounded set $\Om \sub \R^n$. Fix $\varphi \in C^{\infty}_{c}(\Om; \R^N)$ and for sufficiently small $\e _0> 0$ consider the function $\e \mapsto E(u+\e \varphi )$, which is differentiable for $|\e| <\e_0$. We call its derivative $\frac{d}{d\e}E(u_0+\e \varphi )$ at $\e =0$ the first variation of $E$ at $u$ in direction of $\varphi $, and at any critical point $u$ we have $\frac{d}{d\e}\big|_{\e =0}E(u+\e \varphi )=0$.
\[
\underset{\scriptstyle{\text{Figure 2.  }}}{ \includegraphics[scale=.40]{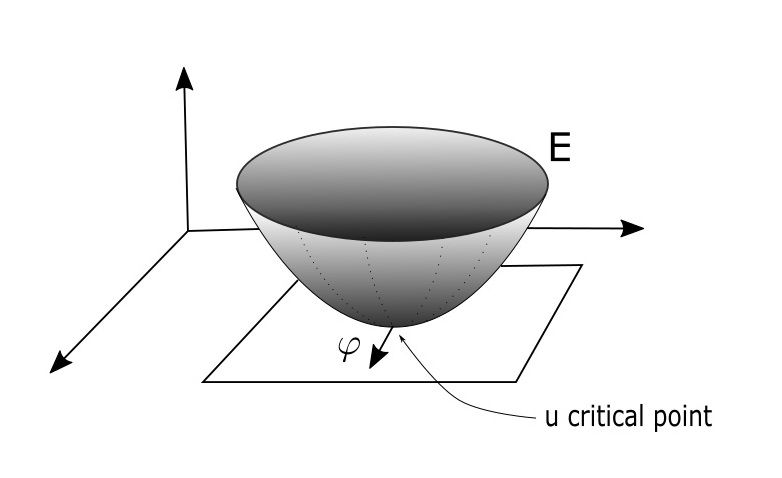}}
\]
Then the directional derivative of $E$ can be computed as follows
 \begin{align}
\frac{d}{d\e}\Big|_{\e =0}E(u+\e \varphi )& =\ \int_\Om \frac{d}{d\e}\Big|_{\e =0} H \big(\cdot,u+\e \varphi,D u+\e D \varphi \big)\   \\  \nonumber
&=\  \int_\Om \Big( H_\eta(\cdot,u,D u) \cdot \varphi + H_\P(\cdot,u,D u) : \D\varphi \Big).
 \end{align}
Here $ H_\eta$ and $ H_\P$ denotes respectively the derivatives of the coefficients $\D_\eta H = (H_{\eta_\al} )_{\al=1,...,N}$ and $\D_\P H = (H_{\P_{\al i}} )_{i=1,...,n}^{\al=1,...,N}$. If $u$ is a global minimiser, namely when
\[
E (u)= \underset{v \in W^{1,p}(\Om)}{\min} E(v),
\]
then $E(u) \leq E(u+\e \varphi)$, for all $\varphi \in C^{\infty}_{c}(\Om; \R^N)$. Therefore,  we have
 \[
\int_\Om \bigg(\sum_{\al , i} H_{\P_{\al i}}(\cdot,u,D u) D_i\varphi _\al +  \sum_{\al} H_{\eta_\al}(\cdot,u,D u) \varphi_\al \bigg)=0, \ \ \ \forall\, \varphi \in C^{\infty}_{c}(\Om; \R^N).
 \]
This is the weak formulation of the Euler-Lagrange PDEs:
\[
\sum_{i} \D_i \big( H_{\P_{\al i}}(\cdot, u, D u) \big) + H_{\eta_\al} (\cdot, u , D u)=0, \ \ \ \al =1,...,N , \ \ \text{in} \ \Om.
\]
This system can be written in compact vector notation as 
  \[
\Div \big( H_\P (\cdot,u,D u)\big)=  H_\eta(\cdot,u,D u), \ \ \ \text{in} \ \Om. 
\]

\subsection{Lower semi-continuity, coercivity and weak compactness }\label{subsection2.4} Non-linearity and weak convergence does not play well together in functional analysis, and this is a source of grave difficulties. For the minimisation problem to be solvable, we need the functional $E: \X \longrightarrow \R$ to satisfy two main properties:
\begin{enumerate}

\item $E$ has to be weakly sequentially lower semi-continuous (LSC) in $\X$, i.e.
\[
E(u) \, \leq\, \underset{m\ri\infty}{\liminf}\, E(u_m), \ \ \text{ when }u_m \weak u  \ \ (\text{or}  \ \ u_m \weakstar \hspace{1pt} u) \text{ as }m\to\infty,
\] 
with respect to the weak topology of the space $\X$ which typically is a reflexive (or a dual) Banach space.

\item We need $E$ to be coercive, namely $E(u) \larrow \infty$ when $\|u \| \to \infty$.
\end{enumerate}
Under these assumptions, we typically have the sequential weak (or weak*) compactness for minimising sequences in $\X$:
\[
E(u_m) \longrightarrow \inf_\X E \ \ \  \Longrightarrow \ \ \  u_{m_k} \weak \hspace{1pt} u \ \ (\text{or}  \ \ u_{m_k} \weakstar \hspace{1pt} u), \text{ as }k\to\infty.
\]
This is sufficient to solve the problem, and comprises the direct method of the Calculus of Variations \cite{D, Ev, KV}, allowing fairly easily to obtain existence of global minimisers. But how one can confirm weak sequential lower semi-continuity?

\subsection{Convexity notions and minimisers}\label{subsection2.5} Consider again $\X= W^{1,p}(\Om; \R^N)$, where $\Om \sub \R^n$ is open and $p \in (1,\infty)$, and the problem of minimising $E:  W^{1,p}(\Om; \R^N) \longrightarrow \R$ given by
\[
E(u)=\int_\Om H(D u)\, dx.
\]
Then, $E$ is sequentially weakly LSC on $W^{1,p}$ if and only if $H$ is (Morrey) quasiconvex on $\R^{N\by n}$. When $\min\{n,N\}=1$ this means that $H$ is convex. If $\min\{n,N\}\geq2$ then it means
\[
\av_\Om H \big( P+D \varphi(x) \big)\, dx  \geq  H(P), \ \ \ \forall \ \varphi \in C^\infty_c(\Om; \R^N).
\]
It is equivalent to requiring that every affine function is a minimiser with respect to its own boundary conditions, and it is non-local condition \cite{D, Ev}. For a general discussion on nonlinear weakly (lower semi-) continuous functionals on Banach spaces we refer e.g.\ to \cite{KV}.

\subsection{Global minimisers are ``absolute" }\label{subsection2.6} A key property of integral functionals is that global minimisers are ``absolute", namely minimise on any subdomain with respect to their own boundary values. This owes to the additivity property of integral. To see this, consider again the integral functional $E$, now defined also on open (or even measurable) subsets:
\[
E(u,\mathcal{O}) :=\int_\mathcal{O} H(\cdot,u,D u)\, dx, \ \ \ \mathcal{O}  \Subset \Om.
\]
Here $\mathcal{O}$ is compactly contained in $\Om$, and $u \in W^{1,p}(\Om; \R^N)$. We will say that \emph{$u$ is an absolute minimiser of $E$ on $\Om$}, if 
\[
E (u,\mathcal{O})\,\leq\, E (u+\phi,\mathcal{O}),\ \ \forall \, \mathcal{O}\Subset \Om, \ \ \forall\, \phi\in W^{1,p}_0(\mathcal{O};\R^N).
\]
Since $W^{1,p}_0(\mathcal{O};\R^N) \sub W^{1,p}_c(\Om;\R^N)$,  any global minimiser $u$ satisfies
\[
\int_\Om H(\cdot,u,D u) \leq \int_\Om H(\cdot,u+\phi,D u + D \phi),
\]
for any $\phi\in W^{1,p}_0(\mathcal{O};\R^N) $. By the additivity property of integrals and minimality, we estimate
 \begin{align}
 \int_{\Om \setminus \mathcal{O}}  H(\cdot,u,D u) & + \int_\mathcal{O} H(\cdot,u,D u)  \nonumber 
\\
& \leq \int_{\Om \setminus \mathcal{O}} H(\cdot,u+\phi,D u + D \phi)  + \int_\mathcal{O} H(\cdot,u+\phi,D u + D \phi)  \nonumber
\\
& = \int_{\Om \setminus \mathcal{O}} H(\cdot,u,D u)  + \int_\mathcal{O} H(\cdot,u+\phi,D u + D \phi) , \nonumber
\end{align}
where in the last step we used that $\phi \equiv 0$ on $\Om\set\mO$. Hence,
\[
\int_\mathcal{O} H(\cdot,u,D u) \leq \int_\mathcal{O} H(\cdot,u+\phi,D u + D \phi) ,
\]
namely $E (u,\mathcal{O})\leq E (u+\phi,\mathcal{O})$, for all $\mathcal{O}\Subset \Om$ and all $\phi\in W^{1,p}_0(\mathcal{O};\R^N)$.


\section{Vectorial Calculus of Variations in $L^\infty$} \label{section3}

\subsection{Supremal functionals and absolute minimisers} A general first order supremal functional has the form
\beq \label{3.1}
E_{\infty}(u,\mO) :=  \underset{x\in\mO}{\ess \, \sup}\ H
\big(x,u(x),Du(x)\big), \ \ \mO  \Subset \Om,
 \eeq
defined on locally Lipschitz maps  $u:\R^n \supseteq \Om\longrightarrow \R^N$, or equivalently for $u \in W_{loc}^{1, \infty}(\Om; \R^N)$. Here $H \in C^1(\R^n \by \R^N \by \R^{N \by n})$ is a (usually non-negative) function called Hamiltonian or supremand, $\Om$ is an open subset of $ \R^n$ and $n,N\in\N$. The simpler and archetypal example of a supremal functional is the ``$\infty$-Dirichlet" functional:
\beq \label{3.2}
\ \ \ E_\infty(u,\mathcal{O}) :=  \| D u\|_{L^\infty(\mathcal{O})}, \ \ \ u \in W^{1,\infty}(\Om;\R^N), \ \mathcal{O} \Subset \Omega,
\eeq
where the $L^\infty$ norm of $Du$ is interpreted as the essential supremum of the Euclidean norm of the gradient matrix. Unlike integral functionals, the lack of additivity for supremal functionals implies that global minimisers are not always absolute! For example consider
\[
E_\infty(u,\mathcal{O})\, :=\, \underset{x\in \mathcal{O}}{\text{ess}\,\text{sup}}\, |u'(x)|^2,\ \ \mathcal{O}\sub \Om,\ \ u \in W^{1,\infty}(\Om),
\]
where $\Om\, =\, (-1,0)\cup(0,1)$. Then there exist two minimisers $u^\pm$ on $\Om$ with boundary data $u(-1)=u(0)=0$, $u(1)=1$ which do not minimise on $(-1,0)$ (Figure 3). Instead, only $u_0\equiv 0$ minimisers thereon. 
\[
\underset{\scriptstyle{\text{Figure 3.  }}}{ \includegraphics[scale=.39]{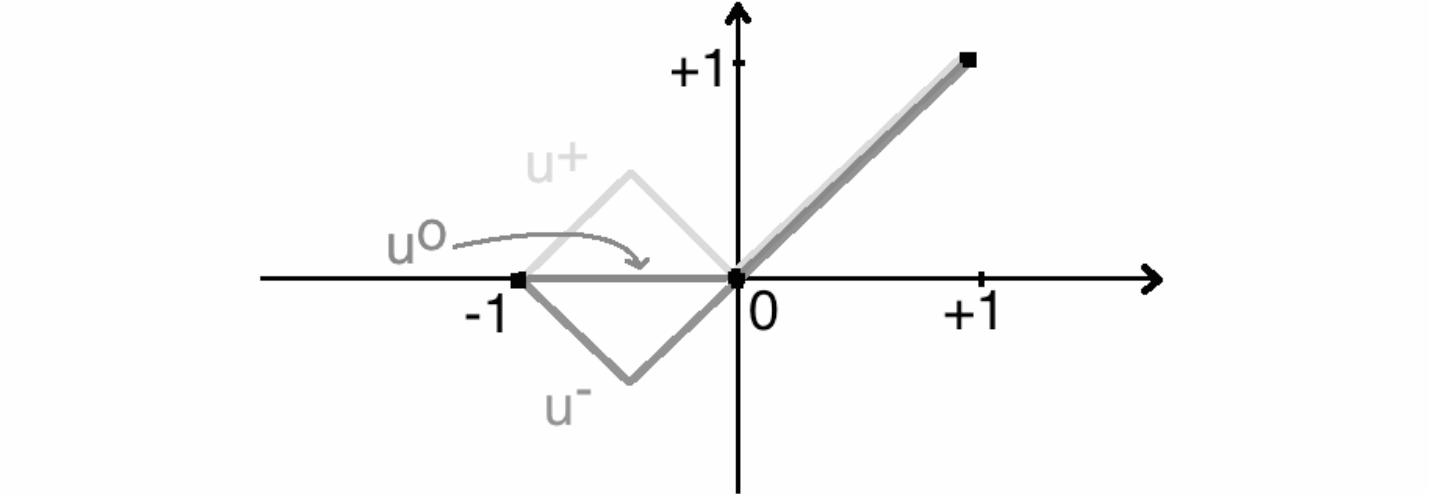}}
\]
Given this lack of locality, the trademark minimality notion in $L^\infty$ is that of \emph{absolute Minimisers}, in which locality is at the outset built into the minimality notion:
\beq \label{1.5}
E_\infty (u,\mathcal{O})\,\leq\, E_\infty (u+\phi,\mathcal{O}),\ \ \forall \, \mathcal{O}\Subset \Om, \ \ \forall\, \phi\in W^{1,\infty}_0(\mathcal{O};\R^N).
\eeq
However, this minimality notion creates more problems than it actually solves, as it throws out of the window the direct method of the Calculus of Variations! Unfortunately, the direct method can only yield global, not absolute minimisers. Further, convexity notions in $L^\infty$ are not as momentous as in the integral case, as weak* LSC plus coercivity does not suffice to provide us with absolute minimisers (see e.g.\ \cite{P, P2, RZ2} for convexity notions in $L^\infty$ and application of weak* lower semi-continuity and $\Gamma$-convergence).

The next natural question would be whether one can replicate the approach of the previous section to derive a PDE satisfied by (absolute) minimisers. Sadly, {\it this is not the case}. The main issue that $E_\infty$ is  \textbf{not} differentiable! This is not a matter related to the properties of $H$ as not even \eqref{3.2} is directionally differentiable for $n=N=1$. The issue stems from the nonlinearity of the ``ess\,sup" operation. (One can confirm this by examples, but this would require a digression to semi-derivatives and Danskin's theorem.) So, if we cannot differentiate $E_\infty$, how could we discover a PDE for (absolute) minimisers, if at all possible?

\subsection{The single Aronsson equation in $L^\infty$}\label{subsection3.1} Let us consider \eqref{3.1} for $N=1$. The idea, which is traced back to Aronsson \cite{A1}-\cite{A7} and is still a trademark feature of the theory nowadays, is to derive a PDE by looking at the limit of the Euler-Lagrange equation of the $p$-Dirichlet functional, or equivalently of the $L^p$-norm of the gradient $\|D \cdot \|_{L^p(\Om)}$. For a fixed function we have $\|Du\|_{L^p(\Om)}\ri \|Du\|_{L^\infty(\Om)}$ as $p\ri \infty$, so there is the expectation that at least at a formal level $\De_p u\ri \De_\infty u$ for some limiting PDE. However it was not a priori clear that the following rectangle ``commutes" 
\begin{align} \label{1.3}
&\|Du\|_{L^p(\Om)}   \ \ \ \  \iff \ \ \ \ \De_p u =0  \nonumber\\
& \ \ \ \downarrow\ p \ri \infty\ \ \ \ \circlearrowleft \ \ \  \ \ \ \ \ \ \downarrow\ p \ri \infty  \\
&\|Du\|_{L^\infty(\Om)}   \ \ \ \overset{?}{\dashleftarrow \dashrightarrow}  \ \ \De_\infty u=0,   \nonumber
\end{align}
and the PDE that one can discover is indeed related to the supremal functional in any way. Remarkably, however, this turned out to be the case! In the limit one can formally discover (the derivation is performed subsequently in the full vectorial case) 
\beq \label{a}
 \De_\infty u:=Du \ot Du :D^2 u=0,
\eeq 
called the \emph{$\infty$-Laplacian}, which is the formal limit of the $p$-Laplacian as $p\to\infty$. Aronsson did indeed prove that, at least in the $C^2$ case, $\infty$-harmonic functions, namely solutions to the $\infty$-Laplacian, can be characterised as the absolute minimisers of the $\infty$-Dirichlet functional  $\|D \cdot \|_{L^\infty(\Om)}$. This result was subsequently extended several years later to the general case of $W^{1,\infty}$ absolute minimisers, by utilising the theory of viscosity solutions \cite{BEJ, C, K7}. The $\infty$-Laplacian is a non-divergence PDE, and weak solutions based on duality methods and integration by parts do not work, whilst on the other hand one can exhibit absolute minimisers which are not $C^2$ differentiable, and for which the PDE makes no rigorous sense. This is one such, known as Aronsson's example: the  function $(x,y)\mapsto|x|^{\frac{4}{3}}-|y|^{\frac{4}{3}}$ is an absolute minimiser, but it cannot be interpreted as a classical solution of $ \De_\infty u=0$, as it is a singular saddle which is not twice differentiable on the axes (Figure 4). 
\[
\underset{\scriptstyle{\text{Figure 4.  }}}{ \includegraphics[scale=.3]{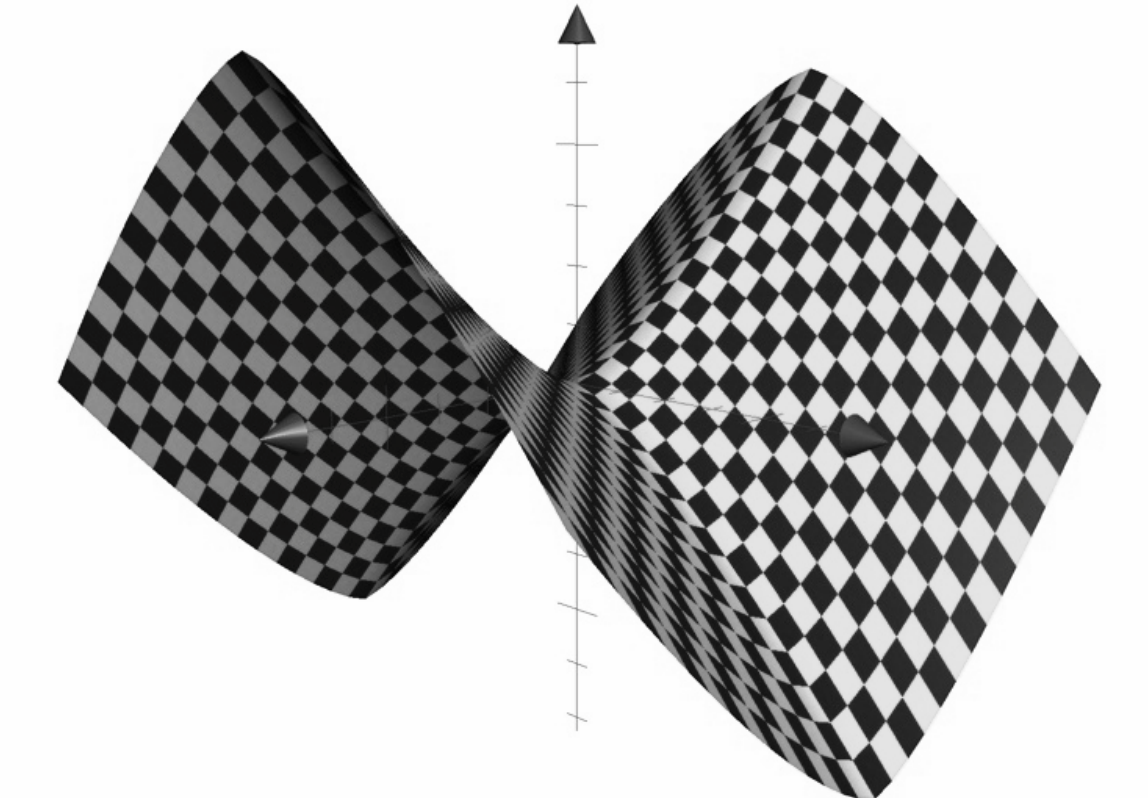}}
\]
More generally, one can apply the above to the general supremal functional \eqref{3.1}, and derive what we now call \emph{Aronsson's PDE} (see e.g.\ \cite{BJW1, BJW2})
 \beq  \label{eq2}
\textrm{A}_\infty u:= H_\P(\cdot,u,D u) \cdot D\big(H(\cdot,u,D u)\big)= 0, \ \ \text{ in }\Om,
 \eeq
which is a second order quasilinear non-divergence degenerate elliptic PDE. Here $H=H(x,\eta,\P)$, and subscripts of $H$ denote derivatives with respect to the corresponding variable. Even though it can be proved that in general absolute minimisers are viscosity solutions to the expanded version of \eqref{eq2}, the converse is not in general true unless very stringent assumptions are satisfied (no $u$-dependence, convexity and uniform coercivity, see \cite{Yu}).

\subsection{A scalar motivational problem in $L^\infty$} (Optimising Lipschitz extensions) Let $\Om \sub \R^n$ be a bounded domain. Suppose $g \in \Lip(\partial \Om)$, and set $L:=\Lip(g, \partial \Om)$, where $\textrm{Lip}$ is the Lipschitz functional
 \beq \label{eq4}
\textrm{Lip}(u,K) := \underset{x,y \in K ,x \neq y}{\sup}
\frac{|u(x)-u(y)|}{|x-y|},\ \ \ K\sub \R^n.
 \eeq
The Lipschitz extension problem asks to find a Lipschitz continuous extension of $g$ with the smallest possible Lipschitz constant, namely to find $u \in \Lip(\Om)$ such that 
\beq \label{1.6}
\left\{
\begin{array}{ll}
 \Lip(u,{\Om})=L,\\
 u=g, \ \ \text{on} \  \partial \Om.
\end{array}
\right.
\eeq
There are multiple solutions to the Lipschitz extension problem, the maximal and minimal of which are the McShane-Whitney extensions (Figure 5):
\beq \label{1.7}
\left\{
\begin{array}{ll}
u^+(x) := \underset {y \in \partial \Om}{\inf} \{ g(y)+L|x-y | \} ,\\
{u}^-(x) := \underset {y \in \partial \Om}{\sup} \{ g(y)-L|x-y | \}.
\end{array}
\right.
\eeq
\[
\underset{\scriptstyle{\text{Figure 5.  }}}{ \includegraphics[scale=.4]{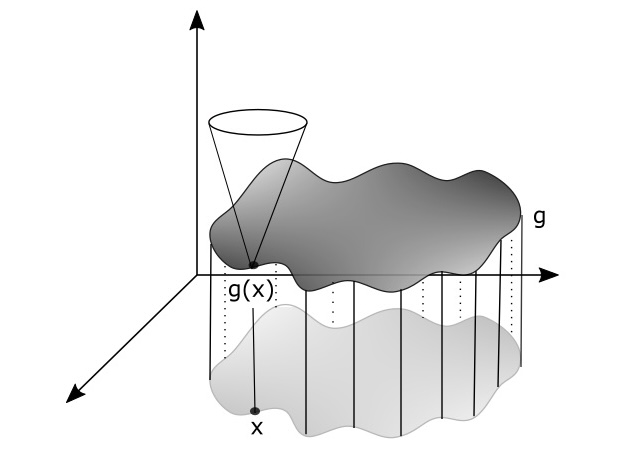}}
\]
Even though it is not a priori obvious, actually one can select a \emph{best} solution to this problem, which is the unique $\infty$-Harmonic function with values equal to those of $g$ (this was Aronsson motivation, see \cite{A3}). The connection between the two problems stems from the identity
\beq \label{Lip}
 \Lip(u,{\Om}) = \| Du\|_{L^\infty(\Om)},
\eeq
which is true for convex domains $\Om$ (or for arbitrary open domains, as long as the Euclidean distance is replaced by the geodesic distance with respect to $\Om$).

\subsection{Aronsson systems in the vectorial case} \label{subsection3.2} Let us now consider the general case of \eqref{3.1} for $N\geq 2$, and mappings $u : \R^n \supseteq \Om \larrow \R^N$. This case is considerably more intricate, and most scalar results and techniques fail in the vectorial case. Let us focus first on the model case of the $\infty$-Dirichlet functional. The method of discovering a PDE by $L^p$ approximations does work, but it reveals more than it does in the scalar case. This is not completely obvious as it needs a renormalisation argument, and before the discovery of the complete $\infty$-Laplace system in \cite{K1}, it was for some time believed that the full $\infty$-Laplacian is only the terms known from the scalar case (see e.g.\ \cite{WO}). To derive (formally) the full system, consider the $p$-Laplacian $\mathrm{div}(|Du|^{p-2}Du)=0$. After distributing the derivatives and normalising, we have the system
\[
Du \ot Du : D^2 u \, +\, \frac{|Du|^2}{p-2}\De u\ = \ 0. 
\]
For any matrix $Q\in \R^{N\by n}$, let $[\![Q ]\!]^\parallel,[\![Q]\!]^\bot$ denote the orthogonal projections on  the range of $Q$ and on its orthogonal complement, when it is seen as a linear map  $Q:\R^n \longrightarrow \R^N$:
\[
\left\{ \ \ \ \ 
\begin{split}
[\![Q ]\!]^\parallel \, &:= \textrm{Proj}_{ \textrm{R}(Q : \R^n \to \R^N)},
\\
[\![Q ]\!]^\bot  &:= \textrm{Proj}_{ \textrm{R}(Q : \R^n \to \R^N)^\bot}.
\end{split}
\right.
\] 
Since $[\![Du]\!]^\parallel + [\![Du]\!]^\bot =\mathrm I$ (the identity on $\R^N$), by expanding the Laplacian $\De u$ with respect to these projections and contracting the derivatives in the first term, we get (see Figure 6)
\beq\label{orth}
Du D\Big(\frac{1}{2}|Du|^2\Big)   \, +\, \frac{|Du|^2}{p-2}[\![Du]\!]^\parallel \De u\ = \ -\frac{|Du|^2}{p-2}[\![Du]\!]^\bot \De u. 
\eeq
\[
\underset{ \ \ \ \text{Figure 6 (a) \hspace{170pt} Figure 6 (b).}}{\includegraphics[scale=0.31]{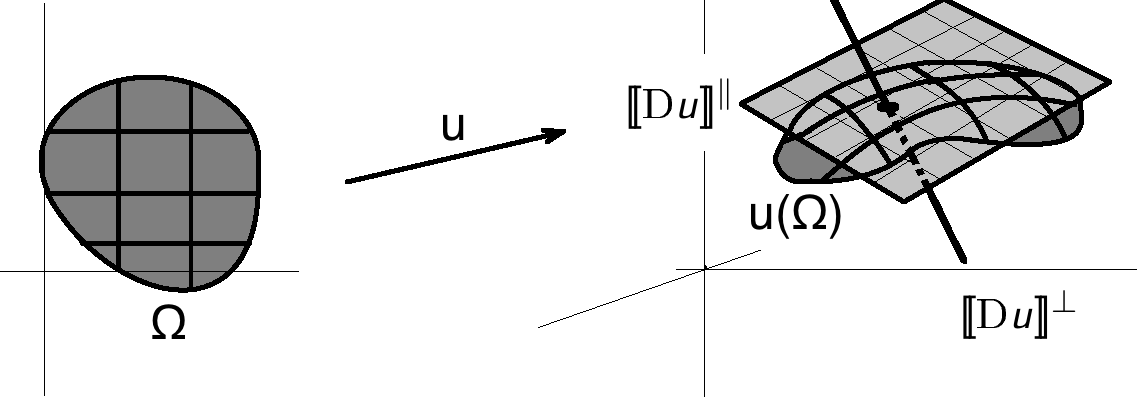}  \includegraphics[scale=0.22]{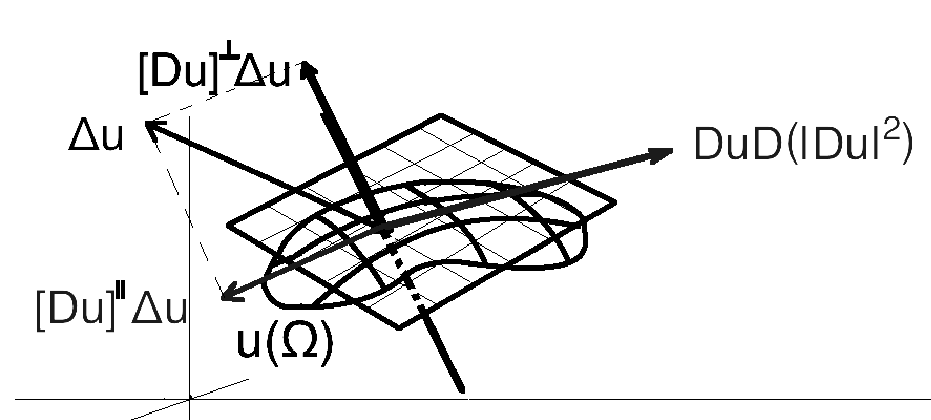}}
\]
By orthogonality, left and right hand sides of \eqref{orth} are orthogonal to each other. Thus, they both vanish and \eqref{orth} decouples to two linearly independent systems. By renormalising the right hand side of \eqref{orth}  and rearranging, we obtain
\[
Du \ot Du : D^2 u \,+\, |Du|^2[\![Du]\!]^\bot \De u\ = \ -\frac{|Du|^2}{p-2}[\![Du]\!]^\parallel \De u. 
\]
Hence, by letting $p\to\infty$ we have derived the complete $\infty$-Laplace system, which reads
\beq 
\label{3.11}
\ \ \Delta_\infty u  := \Big(D u \ot D u + |D u|^2[\![D u]\!]^\bot \! \ot \mathrm{I}\Big):D^{2}u = 0 \ \ \text{ in } \Om.
\eeq
In index from, \eqref{3.11} can be written as
\[
\ \ \ \sum_{\be=1}^N\sum_{i,j=1}^n \Big(\D_i u_\al \, \D_j u_\be + |D u|^2[\![D u]\!]_{\al\be}^\bot\, \de_{ij}\Big)D^2_{ij}u_\be\,=\,0, \  \ \ \al=1,...,N.
\]
A crucial feature of the $\infty$-Laplace system is that it actually consists of two independent systems, one tangential and one normal to the image $u(\Om)$:
\beq \label{split}
D u \, \D\Big(\frac{1}{2}|D u|^2\Big)=0, \ \ \   |D u|^2[\![D u]\!]^\bot \De u = 0.
\eeq
The same approach of $L^p$ approximations can be applied also to the general case of \eqref{3.1}. Suppose first that $H$ depends only on the gradient. Then, the Aronsson system is
 \begin{align}
 \label{3.12}
 \A_\infty u := \Big(H_\P(D u) \ot H_\P(D u) + H(D u)[\![H_\P(D u)]\!]^\bot H_{\P \P}(D u) \Big) :D^2u =  0.
 \end{align}
In the general case that we have dependence of all arguments, the full Aronsson system is
\beq \label{3.13}
\mathcal{A}_\infty \big(\cdot,u,D u,D^2u\big)\,=\,0, \quad \text{in }\Om,
\eeq
where the coefficient map $\mathcal{A}_\infty$ is given by
\beq \label{3.14}
\begin{split}
\mathcal{A}_\infty (x,\eta,P,X) \,  :=&  \Big[  H_P(x,\eta,P) \ot H_P(x,\eta,P)\\
 &  +  \, H(x,\eta,P) [\![H_P(x,\eta,P)]\!]^\bot   H_{PP}(x,\eta,P)\Big]: X  
\\ 
 &+\,  \Big( H_\eta(x,\eta,P) \cdot P \, +\, H_x(x,\eta,P) \Big) H_P  (x,\eta,P)\\
  & + \,  H  (x,\eta,P)[\![H_P(x,\eta,P)]\!]^\bot \Big(H_{P\eta}(x,\eta,P)P
 \\
 & +  \,  H_{Px}  (x,\eta,P) \,-\, H_\eta(x,\eta,P) \Big).
\end{split}
\eeq
Note that the above is just a formal derivation from $L^p$-approximations, and no rigorous direct connection has been established between the Aronsson PDE system \eqref{3.13}-\eqref{3.14} and the supremal functional \eqref{3.1}. Before drawing this connection or discussion existence and other properties, we demonstrate some unexpected properties of   \eqref{3.13}-\eqref{3.14}. Finally, we note that \eqref{3.13}-\eqref{3.14} generally is \emph{not a degenerate elliptic system} for arbitrary $H$, as the coefficient $[\![H_P]\!]^\bot   H_{PP}$ is not in general symmetric, let alone rank-one non-negative (the $\infty$-Laplacian \eqref{3.11} is degenerate elliptic, however). In the scalar case, the Aronsson equation is automatically degenerate elliptic. In \cite{K3} the problem of deciding when the Aronsson system is degenerate elliptic was studied in some detail.


\subsection{Interfaces of discontinuity and phase separation}\label{subsection3.5}

A new difficultly which is not present in the scalar case is that \emph{\eqref{3.13}-\eqref{3.14} has discontinuous coefficients, and this might occur even for smooth solutions}. This happens because there exist solutions whose rank of the gradient is not constant throughout the domain. For simplicity we will limit the discussion to \eqref{3.11}. Large classes of explicit examples are discussed in the next subsection, but here is the simplest one: let $u :\R^2 \larrow \R^2$ be given in complex coordinated by $u(x,y) := e^{ix}-e^{iy}$. Then $u$ is $C^\infty$ and an $\infty$-Harmonic near the origin, having $\rk(Du)=1$ on the diagonal, but $\rk(Du)=2$ otherwise (on $\Om_2$), therefore the projection $[\![D u]\!]^\bot$ is discontinuous (Figure 7). 
\[
\underset{\scriptstyle{\text{Figure 7.  }}}{ \includegraphics[scale=.35]{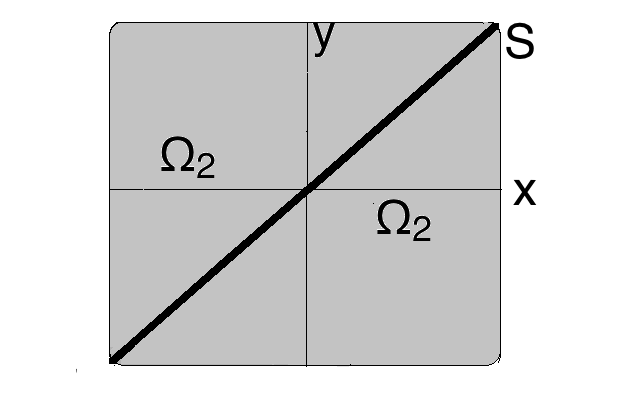}}
\]
Thus, the domain $\Om$ decomposes into two phases $\Om_2$ and $\Om_1$, where $\Om_2\sub \{ \rk(D u) =2 \}$ and $\Om_1 \sub \{ \rk(D u)=1 \}$. This is a general phenomenon, and \emph{$\infty$-Harmonic maps typically present a phase separation}. On each phase the dimension of the tangent space is constant and these phases are separated by \emph{interfaces} whereon the rank of $Du$ ``jumps'' and $[\![Du]\!]^\bot$ becomes discontinuous, whist different behaviours are exhibited on each phase. This phenomenon is well-understood at least in the case of $C^2$ $\infty$-Harmonic mappings $u : \R^2 \supseteq \Om \larrow \R^N$:

\begin{theorem}[Structure of $2$D $\infty$-Harmonic  maps, cf.\ \cite{K3}] Let $ u : \R^2 \supseteq  \Om \longrightarrow \R^{N}$  be an $\infty$-Harmonic  map in  $ C^{2} ( \Om;\R^N )$. Let also $N \geq 2$. Then, there exist disjoint open sets $\Om_{1}$, $\Om_{2} \sub \Om$, and a closed nowhere dense set $\mathcal S$ such that $\Om= \Om_{1} \bigcup \mathcal S \bigcup \Om_{2}$ and:

\smallskip

\noi (i) On $\Om_{2}$ we have  $\rk(D u)=2$, and $ u|_{\Om_{2}}$ is an immersion solution to the Eikonal equation
\[
|D u|^{2}= C \ \ \text{ in }\Om_2,
\]
for some $C\geq0$. The constant $C$ may vary on connected components of $\Om_{2}$.

\smallskip

\noi (ii)  On $\Om_{1}$ we have  $\rk(D u)=1$, and $ u|_{\Om_{1}}$ can be represented by an essentially scalar $\infty$-Harmonic function $f \in C^2(\Om_{1})$, namely
\beq \label{1.11}
u = a+\xi f \ \text{ and } \ \Delta_\infty f=0 \ \ \text{ in }\Om_1,
\eeq
where $a,\xi\in \R^N$. The vectors $a,\xi$ may vary on connected components of $\Om_{1}$.

\smallskip

\noi (iii)  On $ \mathcal S$, $|D u|^{2}$ is constant and also $\rk(D u)=1$. Moreover, if $ \mathcal  S \sub  \partial \Om_1 \cap \partial \Om_2$ (namely if phases coexist) then  $ u|_{\mathcal S}$ is given by an essentially scalar solution to the Eikonal equation.
\end{theorem}

Interestingly, discontinuous coefficients is a genuinely vectorial phenomenon which does not arise when $\min\{n,N\}=1$. Specifically, if $n=1$ then all $\infty$-Harmonic curves are \emph{affine}, and for $u :\R\supseteq \Om \larrow \R^N$, the coefficients in \eqref{3.11} match and it reduces to 
\[
\De_\infty u  = (u' \ot u') u''\, +\, |u'|^2\bigg(I - \frac{u'}{|u'|}\ot \frac{u'}{|u'|}\bigg)u''' = |u'|^2u''.
\]
On the other hand, if $N=1$, then the orthogonal $\infty$-Laplace system $ |D u|^2[\![D u]\!]^\bot \De u =0$ vanishes identically and we reduce to $D u \ot D u :D^2 u=0$.

\subsection{Explicit examples of $\infty$-Harmonic mappings} Both Aronsson's scalar example $u(x,y) = |x|^{4/3} - |y|^{4/3}$ and our earlier vectorial example $u(x,y) = e^{ix}-e^{iy}$ are special case of the following class of solutions $u :\R^2 \supseteq \Om \larrow \R^N$, for $N\geq 1$, having the additively separated form $u(x,y)=f(x) + g(y)$ for  $f,g :\R \larrow \R^N$.  A large family of such separated vectorial solutions was constructed in \cite{K4} and has the form
\beq \label{3.17}
u(x,y):= \int_{y}^{x} e^{iK(t)}\, d t, \ \ \ \  K \in  C^{1}( \R), \ \ \ \| K\|_{L^\infty(\R)} < \frac{\pi}{2}.
\eeq
This class consists of smooth $\infty$-Harmonic maps whose interfaces have triple junctions and even corners. For instance, if $K$ qualitatively behaves as shown in the Figures 8(a), 8(b), then, \eqref{3.17} defines a $C^{2}$  $\infty$-Harmonic  map whose phases are as shown in Figures 9(a), 9(b) respectively. Also, on  $\Om_1$  \eqref{3.17} is given by a scalar $\infty$-Harmonic function, and on  $\Om_2$ it is a solution to the vectorial Eikonal equation. 
 \[
\underset{\scriptstyle{\text{Figure 8 (a). \hspace{100pt} Figure 8 (b).}}}{  \includegraphics[scale=.35]{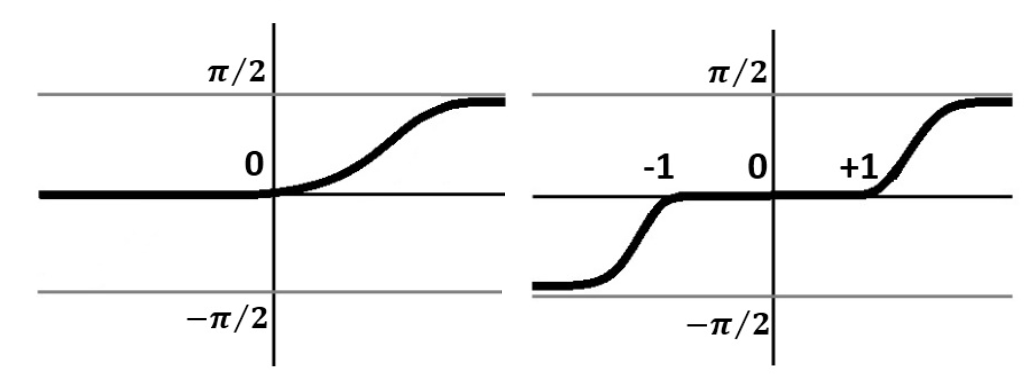} } 
\] 
\[
 \underset{\scriptstyle{\text{Figure 9 (a). \hspace{100pt} Figure 9 (b).\ \ \  }}}{ \includegraphics[scale=.4]{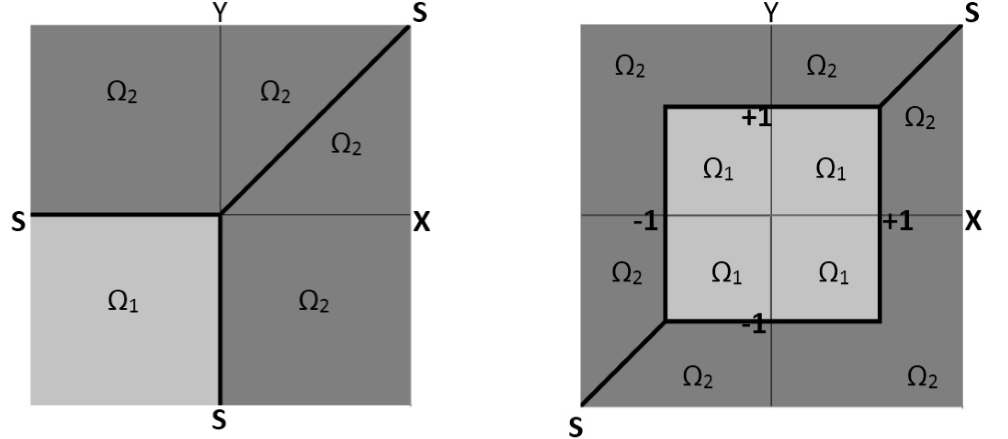}} 
\]
\subsection{Variational structure of $\infty$-Harmonic maps} Recall that the $\infty$-Laplace system \eqref{3.11} has been derived as the limit of $p$-Laplacian as $p\to\infty$, but we have not made any direct connection of \eqref{3.11} to  \eqref{3.2}. In this section we discuss this direct variational connection, taken from \cite{K2}. For the general case of  \eqref{3.1} and \eqref{3.13}-\eqref{3.14}, see \cite{AyK1, AyK2}. Rather surprisingly, unlike the scalar case whereat absolute minimisers (see \eqref{1.5}) is the appropriate minimality notion to characterise variationally $\infty$-Harmonic functions, \emph{if $N\geq2$ vectorial absolute minimisers do not characterise variationally $\infty$-Harmonic maps}. Instead, the appropriate variational notion is different, reflecting the fact that the $\infty$-Laplace system consists of two independent PDE systems, a tangential one and a normal one. Each one of these systems is characterised by a different class of variations: the tangential system is characterised by essentially scalar variations along lines, and the normal system is characterised by variations normal to the image $u(\Om)$ but free on the boundary. Due to technical complications both due to regularity and the discontinuity of the coefficients, this characterisation concerns classical solutions only with constant rank, namely constant rank of the gradient ($\rk(Du)\equiv C$ on $\Om$).

\begin{definition} [$\infty$-Minimal maps cf.\ \cite{K2}] \label{subsection3.7.1}
 Let $u \in C^1(\Om; \R^N)$.

\smallskip

\noi (i) The map $u$ is called a \emph{rank-one absolute minimal} on $\Om$ when for all compactly contained $\mO$ of $\Om$, all $C^1$ functions $f$ on $\mO$ vanishing on $\p \mO$ and all directions $\xi$, $u$ is a minimiser on $\mO$ with respect to essentially scalar variations $u+ \xi f$:
\beq \label{3.6.2}
\left.
\begin{array}{l}
\mO \Subset \Om, \\
f\in C^1_0(\mO), \\
\xi \in \mS^{N-1}
\end{array}
\right\} \ \ \Longrightarrow \ \
E_\infty(u,\mO)\ \leq \ E_\infty(u+f\xi,\mO).
\eeq
\[
\underset{\scriptstyle{\text{Figure 10.  }}}{ \includegraphics[scale=.3]{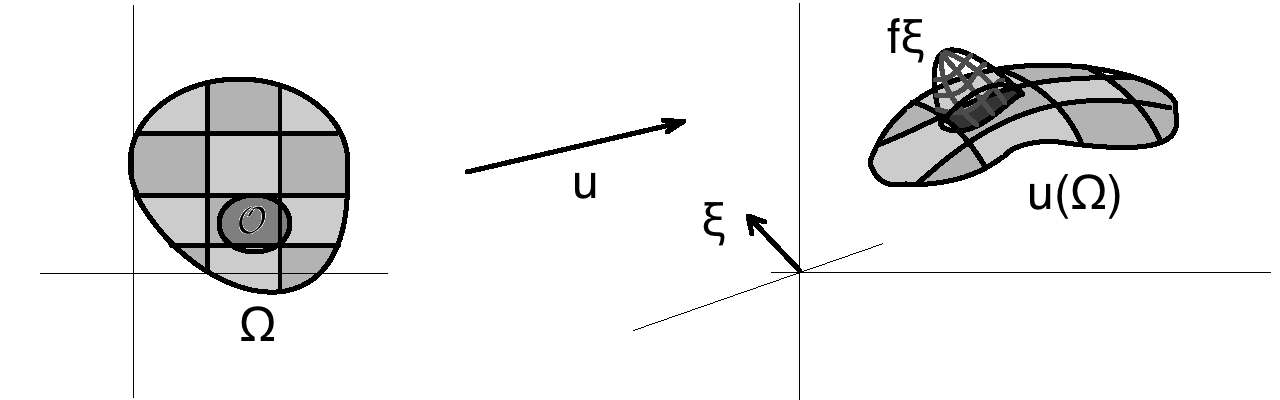}}
\]
 \noi (ii) Suppose that $u$ has constant rank. We say that \emph{$u(\Om)$ has $\infty$-minimal area} when for all compactly contained  $\mO$, all $C^1$ functions $h$ on $\bar{\mO}$ (not vanishing on $\p \mO$) and all normal vector fields $\nu$, $u$ is a minimiser on $\mO$ with respect to normal free variations $u+h\nu$:
\beq \label{3.6.3}
\left.
\begin{array}{l}
\mO \Subset \Om \set \S, \\
h\in C^1(\bar{\mO}), \\
\nu \in \Gamma([\![D u]\!]^\bot)
\end{array}
\right\} \ \ \Longrightarrow \ \
E_\infty(u,\mO)\ \leq \ E_\infty(u+h \nu,\mO).
\eeq
\[
\underset{\scriptstyle{\text{Figure 11.  }}}{ \includegraphics[scale=.3]{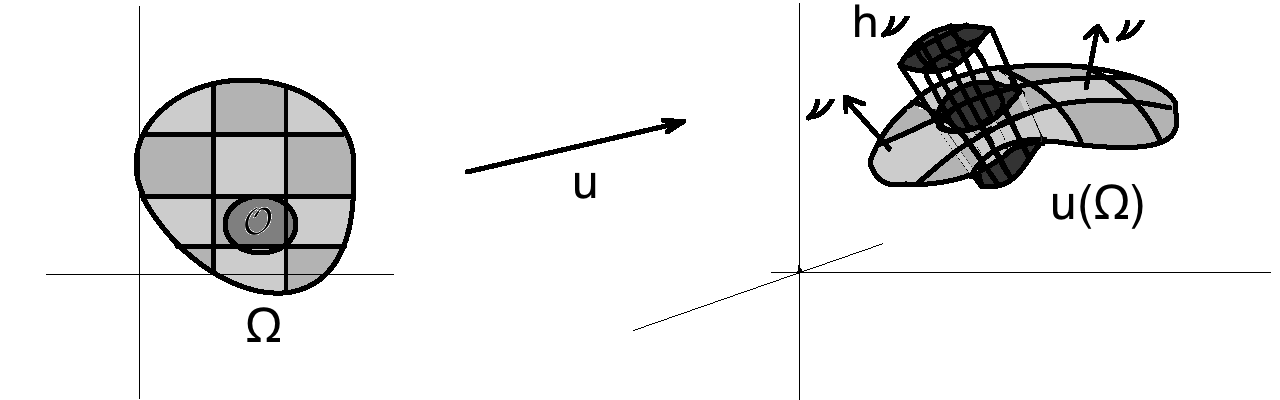}}
\]
\noi (iii) If $u$ is a constant-rank map on $\Om$, we call $u$ an \emph{$\infty$-Minimal map with respect to functional \eqref{3.2}} when $u$ is a rank-one absolute minimal on $\Om$ and $u(\Om)$ has $\infty$-minimal area.
\end{definition}

By utilising the above variational notion, one can establish the following result.

\begin{theorem} [Variational Structure of the $\infty$-Laplacian cf.\ \cite{K2}]  Let $u  \in C^2(\Om;\R^N)$. Then, if $u$ is an $\infty$-Minimal map with respect to \eqref{3.2}, it follows that $u$ is $\infty$-Harmonic and solves \eqref{3.11} on $\Om$. If additionally $u$ is an immersion, the converse is also true and $\infty$-Harmonicity implies $\infty$-Minimality. 

\smallskip 

Additionally, we have the following characterisations for the tangential and the normal components of the $\infty$-Laplace PDE system separately:

\smallskip 

\noi $\bullet$ If $u$ is a Rank-one absolute minimal on $\Om$, then it is tangentially $\infty$-Harmonic  and solves $Du \ot Du :D^2u= 0$ on $\Om$. The converse is true if $u$ is an immersion. 

\smallskip 

\noi  $\bullet$ If $u$ is an immersion, then $u(\Om)$ has $\infty$-Minimal area if and only if it is normally $\infty$-Harmonic  and solves $|Du|^2[\![Du]\!]^\bot \De u= 0$ on $\Om$.
\end{theorem}

We note that, neither the above result nor the extensions in \cite{AyK1,AyK2}, exactly disprove that \eqref{3.11} could perhaps still be deducible from vectorial absolute minimisers, but the \emph{standing conjecture is that this is not possible}. Notwithstanding, the reverse implication is definitely not true: in \cite{KS} \emph{classes of $2$-dimensional smooth $\infty$-Harmonic maps were constructed, which are not even globally minimising for \eqref{3.2}}. Thus, the \emph{$\infty$-Laplace system does not suffice for absolute minimality in $L^\infty$}.

\subsection{Vectorial motivational problems in $L^\infty$} We now discuss two interesting problems for which vectorial $L^\infty$ problems and Aronsson systems are relevant. There are numerous other applications, but most either concern higher problems, or require constraints (see e.g.\ \cite{CKM, CK, GNP, K9}, \cite{K11}-\cite{K14}, \cite{KM3}). It's also worth noting that the scalar motivation of Aronsson on the Lipschitz extension problem does not apply to the vectorial case. Unless one replaces the Euclidean norm on $\R^{N\by n}$ by the operator norm, \eqref{Lip} is no longer true. Using the operator norm as supremand leads to an Aronsson system with non-uniquely defined coefficients, as one needs to select a branch of the eigenvector, see \cite{K1, SS}.

\subsubsection{Optimising quasi-conformal mappings (cf.\ \cite{K5})}  Let $N\geq n\geq 3$. A $C^1$ map  $u :\R^n \supseteq \Om \longrightarrow \R^N$ is conformal if and only if there exists a function $\la$ such that $(D u)^{\top} D u = \la \mathrm{I}_n$, on $\Om$, from where it follows that $\la = \frac{1}{n}|D u|^2$. On the other hand, an immersion $u$ called quasi-conformal if and only if the dilation $K(Du)$ of $u$ is essentially bounded on $\Om$. Here the dilation $K : \R^{N\by n} \larrow [0,\infty]$ is the non-convex function defined as
\[
K(P):=\left\{ \begin{array}{ll}
\dfrac{|P|^2}{\det(P^\top\! P)^{1/n}}, & \mathrm{rk}(P)=n,
\\
+\infty,  &  \mathrm{rk}(P)<n.
\end{array}
\right.
\]
Note that we always have that $K(u) \geq \sqrt{n}$, whilst, by the arithmetic-geometric mean inequality, we have $K(u) = \sqrt{n}$ a.e.\ on $\Om$ if and only if $u$ is conformal. Hence, the dilation measures the deviation of $u$ from being conformal, and the smaller the value of
\beq  \label{3.3.3.1} 
\mathrm K_\infty (u,\Om) := \big\|K(Du)\big\|_{L^\infty(\Om)}
\eeq
is, the closer $u$ is to being a conformal map. However, the gap between conformal and quasi-conformal maps is very big, and in geometric analysis it is of interest to close this gap. Hence, $L^\infty$ vectorial variational methods can be used to optimise quasi-conformal maps \cite{K5}. The Aronsson PDE system associated to \eqref{3.3.3.1} has the form of the quasilinear system \eqref{3.12}, but the exact terms are rather complicated to express. When $n=N$, this approach has previously been followed by Capogna-Raich \cite{CR} and relates to Teichm\"uller's theory. In particular, a conjecture appearing in \cite{CR} was disproved. 

\subsubsection{Variational data assimilation and weather forecasting (cf. \cite{K9, AK})}  An important applications of \eqref{3.1} is to problems of $L^\infty$-modelling of variational Data Assimilation (\emph{4DVar}) arising in the Earth Sciences (see e.g.\ \cite{BS}). An explicit model of $H$ is
\beq \label{3.3.2.1}
H(x,\eta,P) := \big|P-\mathcal{V}(x,\eta) \big|^2 +\, \big|k(x)-K(\eta) \big|^2 ,
\eeq
and describes the ``error" in the following sense: consider the problem of finding the  solution $u$ to the following ODE system coupled by a pointwise constraint
\[
\ \ \ D u(t) = \mathcal{V}\big(t,u(t)\big) \ \ \ \& \ \ \ K(u(t)) = k(t), \ \ \ t\in \Om.
\]
Here  $\mathcal{V} : \Om \by \R^N\larrow \R^N$ is a vector field giving the law of motion along the trajectory described by  $u:\R\supseteq \Om \larrow \R^N$ (e.g.\ Galerkin approximation of the Euler equations, Newtonian forces,  etc), $k: \R \supseteq \Om \larrow \R^M$ is some partial fuzzy (approximate) measurements in continuous time along the trajectory, and $K: \R^N \larrow \R^M$ is a submersion which corresponds to some component of the trajectory that can be observed. The problem is clearly overdetermined as $k$ is known only approximately up to some error. Thus, one can try to minimise the $L^\infty$ error, and minimisation of \eqref{3.1} with $H$ as in \eqref{3.3.2.1} will lead to a uniformly best approximation without ``spikes" of large deviation of the prediction.

\subsection{Motivation of generalised concepts: $\mD$-solutions}\label{subsection3.8} So far we have not discussed the question of existence of solutions to the Dirichlet problem for the Aronsson system \eqref{3.13}-\eqref{3.14}, neither the question of existence of vectorial absolute minimisers with prescribed boundary values. Both problems are highly non-trivial, but since the latter problem is not directly related to the Aronsson system, let us discuss the former question. The latter problem is anyway still open in the vectorial case when $\min\{n,N\}\geq 2$; see \cite{AK, BJW1, C, K7, K9} for the case $\min\{n,N\}=1$ and work towards the general case $\min\{n,N\}\geq 2$ in \cite{KM2}.

The starting difficulty is that \eqref{3.13}-\eqref{3.14} is a non-divergence non-monotone system (with discontinuous coefficients), so none of the classical theories of weak, distributional or viscosity solutions apply. The same observations are true also for the simplest non-trivial model case of the $\infty$-Laplacian \eqref{3.11}. In the scalar case of $N=1$, viscosity solutions are applicable and this is the standard notion of generalised solution for \eqref{a}-\eqref{eq2}. The need for generalised solutions easily becomes apparent. Aronsson's example is not twice differentiable on the axes, and in view of \eqref{split}, any $C^1$ submersion with constant modulus of the gradient formally solves \eqref{3.11}, without actually being able to expand the PDE system.

Motivated by the issues arising in $L^\infty$, a new theory of generalised solutions was proposed in \cite{K8}, which applies to fully nonlinear PDE systems of any order and allows for merely measurable maps as solutions. This approach is not based either on integration by parts or on the maximum principle. Instead, the starting point is a probabilistic representation of derivatives via limits of difference quotients in compactified Young measure spaces.

To explain the underlying ideas, let us motivate the notion in the case of first order ODE systems. The situation will not be too dissimilar to the case of \eqref{3.13}-\eqref{3.14}, as we know that a priori solutions have to be once weakly differentiable in $W^{1,\infty}(\Om;\R^N)$, so we only need to interpret in an appropriate weak sense only the Hessian.

We start with some notation. Given $a\in \R^n$ with $|a|=1$ and $h\in \R\set \{0\}$, consider the difference quotient of a map $u : \R^n \supseteq \Om \larrow \R^N$ by
\beq  \label{3.7.8}
D^{1,h} u(x) := \bigg( \frac{u(x+he_1)-u(x)}{h}, \ldots, \frac{u(x+he_n)-u(x)}{h}\bigg) \in \R^{N\by n}.
\eeq
Let now $n=1$ and consider a fully nonlinear first order ODE system for $u: (a,b) \larrow \R^N$:
\beq\label{3.7.10}
F(\cdot,u,u')=0, \ \  \text{ in } (a,b),
\eeq
where $F: (a,b) \times \R \times \R \larrow \R^N$ is assumed to be continuous. Suppose $u: (a,b) \larrow \R^N$ is a classical (or strong a.e.) solution. Since the derivative is the pointwise (or a.e.) limit of the difference quotients, \eqref{3.7.10} is equivalent to
\beq 
\label{3.7.12}
 F\Big(\cdot,u,\lim_{h\ri \infty} D^{1,h} u\Big)=0, \quad \text{in } (a,b).
\eeq
Since by assumpiton $F$ is continuous, \eqref{3.7.12} is further equivalent to
\beq\label{3.7.11}
\lim_{h\ri \infty} F\big(\cdot,u,D^{1,h} u \big)=0, \quad \text{in } (a,b).
\eeq
The crucial observation is that \eqref{3.7.11} makes sense for any $u: (a,b) \larrow \R^N$, without requiring differentiability, and broadly could be taken as a notion of generalised solution. The motivation itself above shows that for any point that $u$ happens to be differentiable, it satisfies the equation in the classical pointwise sense. The underlying idea is that one replaces the pointwise value of $u'(x)$ (which may not exist), by the asymptotic limit of the difference quotients  (see Figure 12), but inside the coefficient map. 
\[
\underset{\ \ \ \text{Figure 12.}}{\includegraphics[scale=0.35]{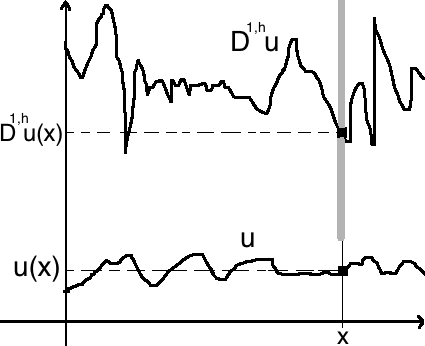}}
\]
Despite being reasonable, the above approach need to be slightly modified, in order to be able to provide good convergence properties for the notion of generalised solution. To motivate how one can represent the limit \eqref{3.7.11} in a useful way, note that \eqref{3.7.10} holds if and only if:
\beq
 \label{3.7.14}
\int_{{\R^N}} \Phi(X)\hspace{1pt}\mF(x,u(x),X)\hspace{1pt} d\big[\de_{(u'(x))} \big](X)\hspace{1pt}=\hspace{1pt} 0, \quad x\in(a,b),
\eeq  
for any compactly supported $\Phi \in C_c\big( \R^N \big)$. Namely, we switch from the classical viewpoint of the derivative as a map $u' : (a,b) \larrow \R^{N}$ by seeing it a probability-valued map given by the Dirac mass at $u'$:
\[
\ \ \de_{u'}\ : \ \ (a,b)  \larrow \mathscr{P}(\R^{N}),\quad    x\mapsto \de_{u'(x)}.
\]
Further, we may restate that $u'$ is the limit (in measure) of the difference quotients $D^{1,h}u$ as $h\ri 0$ by writing
\[
\de_{D^{1,h}u} \weakstar\hspace{1pt} \de_{D u}, \quad \text{as }h\ri0.
\]
The weak* convergence is meant in the space of Young measures valued into $\R^{N}$, that is the set of measurable probability-valued mappings $(a,b) \larrow \mathscr{P}(\R^{N})$. The rationale of the reformulation is that \emph{we may thus allow for general probability-valued ``diffuse gradients" of measurable maps which may not be concentration measures.} This is indeed possible if we replace $\R^{N}$ by its $1$-point compactification $\smash{\overline{\R}}^{N} := \R^{N}\cup\{\infty\}$. So by considering instead the maps $(\de_{D^{1,h}u})_{h\neq 0}$ as Young measures, we obtain the necessary compactness and  we \emph{always} have subsequential weak* limits in the space of Young measures $\mathscr{Y}\big((a,b); \smash{\overline{\R}}^{N}\big)$:
\beq  \label{3.7.6}
\de_{D^{1,h_i}u} \weakstar\hspace{1pt} \mD u, \quad \text{as }h_i\ri0.
\eeq
These will be our generalised derivatives. Then, we may interpret \eqref{3.7.10} for measurable maps $u : (a,b) \larrow \R^N$ as
\beq  \label{3.7.7}
\int_{\smash{\overline{\R}}^{N}} \Phi(X)\hspace{1pt} \mF\big( x,u(x),X\big)\hspace{1pt} d[\mD u(x)](X)\hspace{1pt} =\hspace{1pt} 0, \quad \text{ a.e. }x\in (a,b),
\eeq
 for any test function $\Phi \in C_c\big( \R^{N} \big)$ and any ``diffuse derivative" $\mD u \in \mathscr{Y}\big((a,b); \smash{\overline{\R}}^{N}\big)$.  \eqref{3.7.6} and \eqref{3.7.7} essentially constitute the definition of \textbf{diffuse derivatives} and \textbf{$\mD$-solutions} in the special case of \eqref{3.7.10} and is the central notion of solution for Aronsson systems.

\begin{definition}[Weakly differentiable $\mD$-solutions of $2$nd order PDE systems, cf.\ \cite{K8}]
 Let $\Om \sub \R^n$ be open, $\mF : \Om \by \big(\R^N\by \R^{Nn}\by \R^{Nn^2}_s \big) \larrow \R^M$ a Carath\'eodory map and $u : \R^n \supseteq \Om \larrow \R^N$ a map in $W^{1,1}_{\text{loc}}(\Om;\R^N)$. Consider the PDE system
\beq \label{3.8.1.1}
\mF\big(\cdot,u,D u,D^2 u\big)\hspace{1pt} =\hspace{1pt}0, \ \ \text{ in }\Om.
\eeq
We will say that $u$ is a \textbf{$\mD$-solution of \eqref{3.8.1.1}} when for any diffuse hessian $\mD ^2 u \in \mY(\Om,\smash{\overline{\R}}^{Nn^2}_s)$ of $u$ and any $\Phi \in C_c\big( {\R}^{Nn^2}_s \big)$ we have
\beq \label{3.8.1.2}
\ \ \ \int_{\smash{\overline{\R}}^{Nn^2}_s} \Phi(\X)\hspace{1pt} \mF\big( \cdot,u,D u,\X \big)\hspace{1pt} d[\mD^2 u ](\X)\hspace{1pt} =\hspace{1pt} 0, \ \ \text{ a.e.\ on } \Om.
\eeq
\end{definition}
The general properties of Young measures imply that $u$ is differentiable in measure with derivative $D u$ if and only if the diffuse gradient $\mD u \in \mY(\Om,\smash{\overline{\R}}^{Nn})$ is unique and $\mD u = \de_{D u}$ a.e.\ on $\Om$. For details on the comprehensive theory of $\mD$-solutions and applications, see \cite{K8}.

\begin{remark}[Young measures] We recall some elements of Young measures for the convenience of the reader (see \cite{FG, K8}). The $L^1$-space of maps valued into a separable Banach space $E$ is
\[
L^1\big( \Om, E\big) \hspace{1pt} :=\hspace{1pt} \Big\{ \Phi:\Om \larrow E\hspace{1pt} : \hspace{1pt} \Phi \text{ is strongly measurable and } \| \Phi \|_{L^1( \Om, E)} < \infty \Big\}.
\]
In particular, if $E=C(\mK)$, the elements $\Phi$ of $ L^1\big( E, C(\mK)\big)$ are those Carath\'eodory functions for which $\| \Phi \|_{L^1( E, C(\mK))}:= \int_E \max_{X\in \mK} \big|\Phi(x,X)\big|\hspace{1pt} dx < \infty$. It is well-known (see e.g.\ \cite{FL}) that
\[
\left( L^1\big( E, C(\mK)\big) \right)^* \hspace{1pt} =\hspace{1pt} L^\infty_{w^*}\big( E,\mM(\mK) \big).
\]
The dual space above consists of measure-valued maps $E \ni x   \mapsto \vartheta(x) \in \mM(\mK)$ which are weakly* measurable, i.e.\ for any open $\mathcal{U} \sub \mK$, the function $[\vartheta(\cdot)](\mathcal{U}) \in \R$ is measurable on $E$:
\[
 L^\infty_{w^*}\big( E,\mM(\mK) \big) \hspace{1pt} :=\hspace{1pt} \Big\{\vartheta:E \larrow \mM(\mK) \hspace{1pt} : \hspace{1pt} \vartheta \text{ is weakly* measurable and }  \| \vartheta \|_{L^\infty_{w^*} ( E,\mM(\mK) )} < \infty \Big\}.
\]
The norm is $\| \vartheta \|_{L^\infty_{w^*} ( E,\mM(\mK) )} := {\ess\hspace{1pt}\sup}_{x\in E} \left\|\vartheta(x) \right\|$, where ``$\|\cdot\|$" denotes the total variation. Since $L^1\big( E, C(\mK)\big)$ is separable, the closed unit ball of $L^\infty_{w^*}\big( E,\mM(\mK) \big)$ is sequentially weakly* compact.  The duality pairing $\langle\cdot,\cdot\rangle :  \ L^\infty_{w^*}\big( E,\mM(\mK) \big) \by L^1\big( E, C(\mK)\big) \larrow \R$ is given by
\[
\langle \vartheta, \Phi \rangle\hspace{1pt} :=\hspace{1pt} \int_E \int_{\mK} \Phi(x,X)\hspace{1pt} d[\vartheta(x)] (X)\hspace{1pt} dx.
\]
Let $E\sub \R^n$ be a measurable set and $\mK$ a compact subset of some Euclidean space $\R^d$. The subset of the unit sphere of $L^\infty_{w^*}\big( E,\mM(\mK) \big)$ which consists of probability-valued maps is called the set of Young measures:
\[
\mY(E,\mK)\hspace{1pt} :=\hspace{1pt} \Big\{ \vartheta\hspace{1pt} \in \hspace{1pt} L^\infty_{w^*}\big( E,\mM(\mK) \big)\hspace{1pt} : \hspace{1pt} \vartheta(x) \in \mP(\mK),\text{ for a.e. }x\in E\Big\}.
\]
\end{remark}

Now we close these notes with an existence result for $\mathcal D$-solutions to the Dirichelt problem for the Aronsson system, first established in the special case of the $\infty$-Laplacian in \cite{K8}, and subsequently extended to more general systems in \cite{CKP}.

\begin{theorem}[Existence of $\mD$-solutions to the Dirichlet problem for the Aronsson system, cf.\ \cite{CKP, K8}] \label{theorem3} Let $N\leq n$ and suppose $H : \R^{N\by n} \larrow [0,\infty)$ satisfies $H(P) =h\big(PP^\top\big)$ for some $h : \R^{N\by N}_{\geq 0} \larrow [0,\infty)$ such
that $h \in C^1\big(\R^{N\by N}_{\geq 0}\big) \cap C^2\big(\R^{N\by N}_{> 0}\big)$, the derivative $h_X$ is symmetric, and $\det(h_X(X))\neq 0$ for $X>0$. (These assumptions are automatically satisfied for \eqref{3.11}.) Then, for any open set $\Om\sub \R^n$ and $g\in W^{1,\infty}(\Om,\R^N)$, the Dirichlet problem for \eqref{3.12}
\begin{equation}\label{3.31}
\left\{
\begin{array}{rl}
\mathrm{A}_\infty  u  = 0, & \text{ in }\Om,\\
u= g,  & \text{ on }\p\Om,
 \end{array}
 \right.
\end{equation}
has an infinite set of $\mD$-solutions in the class of Lipschitz submersions 
\[
\mA\,:=\, \Big\{v\in W^{1,\infty}_g(\Om,\R^N)\ :\ \rk(\D v)=N, \text{ a.e.\ on }\Om \Big\}.
\]
Namely, there exists an infinite set of $u\in \mA$ such that, for any $\mD ^2 u \in \mY(\Om,\smash{\overline{\R}}^{N\by n^2}_s)$,
\[
\ \ \ \int_{\smash{\overline{\R}}^{N\by n^2}_s} \Phi(\X)\, \Big(H_P \ot H_P + H[\![H_P]\!]^\bot H_{PP}\Big)(\D u):\X\, d[\mD^2u](\X)\, =\, 0
\] 
a.e.\ on $\Om$, for any $\Phi \in C_c\big( {\R}^{N\by n^2}_s \big)$.
\end{theorem}

Unlike the scalar case (see e.g.\ \cite{C, K7}), the Dirichlet problem \eqref{3.31} has a non-unique solution, but this is not a defect of the notion of generalised solution. In \cite{KS} it is shown that not even smooth $\infty$-Harmonic maps on the unit disc are unique for $n=N=2$. The method of proof for Theorem \ref{theorem3} is based on the use of the Dacorogna-Marcellini Baire category method \cite{DM} for an associated differential inclusion, which can be understood any analytic counterpart of Gromov's geometric Convex Integration method.

\bibliographystyle{amsplain}

\end{document}